\documentclass[11pt,twoside]{article}
\usepackage{amssymb}
\newcommand{\ignore}[1]{}
\ignore{
The text in between these parentheses is ignored by LaTeX
ignore}

\def\@begintheorem#1#2{\par\bgroup{\sc #1\ #2. }\it\ignorespaces}
\def\@opargbegintheorem#1#2#3{\par\bgroup{\sc #1\ #2\ (#3). } \it\ignorespaces}
\def\@endtheorem{\egroup}
\newtheorem{theorem}{Theorem}[section]
\newtheorem{corollary}[theorem]{Corollary}
\newtheorem{lemma}[theorem]{Lemma}
\newtheorem{proposition}[theorem]{Proposition}
\newtheorem{example}[theorem]{Example}
\newtheorem{algorithm}[theorem]{Algorithm}
\newtheorem{definition}[theorem]{Definition}
\newcommand{\bt}[1]{\begin{theorem}\label{#1}}
\newcommand{\bc}[1]{\begin{corollary}\label{#1}}
\newcommand{\bl}[1]{\begin{lemma}\label{#1}}
\newcommand{\bp}[1]{\begin{proposition}\label{#1}}
\newcommand{\be}[1]{\begin{example}\rm\label{#1}}
\newcommand{\ba}[1]{\begin{algorithm}\rm\label{#1}}
\newcommand{\bd}[1]{\begin{definition}\rm\label{#1}}
\newcommand{\bpr}{\par{\it Proof}. \ignorespaces}
\newcommand{\et}{\end{theorem}}
\newcommand{\ec}{\end{corollary}}
\newcommand{\el}{\end{lemma}}
\newcommand{\ep}{\end{proposition}}
\newcommand{\ee}{\end{example}}
\newcommand{\ea}{\end{algorithm}}
\newcommand{\ed}{\end{definition}}
\newcommand{\epr}{{\ \vbox{\hrule\hbox{%
\vrule height1.3ex\hskip0.8ex\vrule}\hrule
}}\\\par}

\def\Q{\mathbb{Q}}
\def \C {{{\cal C}}}
\def \F {{{\cal F}}}
\def \P {{\cal P}}
\def \Z {{\cal Z}}
\def \conv {\mbox{conv}}
\def \zone {\mbox{zone}}

\begin{document}

\title{\bf Convex Combinatorial Optimization}
\author{Shmuel Onn
\thanks{Supported in part by a grant from ISF - the
Israel Science Foundation, by a VPR grant at the Technion,
by the Fund for the Promotion of Research at the Technion,
and at MSRI by NSF grant DMS-9810361.}
\and Uriel G.\ Rothblum
\thanks{Supported in part by a grant from ISF - the
Israel Science Foundation, by a VPR grant at the Technion,
and by the Fund for the Promotion of Research at the Technion.}
}

\date{}
\maketitle

\centerline{\em Dedicated to Professor Louis J.\ Billera on the
Occasion of his Sixtieth Birthday}

\begin{abstract}
We introduce the convex combinatorial optimization problem, a far
reaching generalization of the standard linear combinatorial
optimization problem. We show that it is strongly polynomial time
solvable over any edge-guaranteed family, and discuss several
applications.
\end{abstract}

\section{Introduction}

The general linear combinatorial optimization problem is the following.

\vskip.2cm\noindent{\bf Linear Combinatorial optimization.} Given
a family $\F\subseteq 2^N$ of subsets of $N:=\{1,\dots,n\}$ and a
rational weighting $w:N\longrightarrow\Q$, find $F\in\F$ of
maximum weight $w(F):=\sum_{j\in F}w(j)$.

\vskip.4cm\noindent
There is a massive body of knowledge on the computational complexity
of this problem for various classes of families presented in various ways
(in terms of $n$ and sometimes additional parameters),
and efficient algorithms in numerous cases, cf. \cite{CCPS}.
For instance, if $\F$ is the family of stable sets in a given graph
with vertex set $N$ then the problem is NP-hard
whereas if $\F$ is the family of matchings in a given graph with edge set $N$
then the problem is polynomial time solvable.

\vskip.4cm In this article we consider the following
generalization of linear combinatorial optimization.

\vskip.2cm\noindent{\bf Convex combinatorial optimization.} Given
$\F\subseteq 2^N$ with $N=\{1,\dots,n\}$, a
vectorial weighting $w:N\longrightarrow\Q^d$, and a convex functional
$c:\Q^d\longrightarrow \Q$, find $F\in\F$ of maximum value $c(w(F))$.

\vskip.4cm\noindent
The standard linear combinatorial optimization problem over a family $\F$
is recovered as the special case with $d=1$, $w:N\longrightarrow\Q$
weighting by scalars, and $c:\Q\longrightarrow\Q:x\mapsto x$ the identity.

\vskip.4cm
Convex combinatorial optimization has a very broad expressive power and
conveniently captures a variety of problems studied in the
operations research and mathematical programming literature
including quadratic assignment, inventory management, scheduling,
reliability, bargaining games, clustering, and vector partitioning,
see \cite{AFLS, AAS, BHR, BH, COR, GR, HHPR, PRW} and references therein.
In Section 3 we discuss some of these applications in detail and
demonstrate that, as a consequence of our framework,
all admit a simple unified strongly polynomial time algorithm.

\vskip.4cm A particularly successful general methodology for
linear combinatorial optimization is the geometric
approach inaugurated by Edmonds \cite{Edm1} and culminated in
Gr\"otschel-Lov\'asz-Schrijver \cite{GLS}, \break
outlined as follows.
With each family member $F\in\F$ is associated its indicator
${\bf 1}_F:=\sum_{j\in F} {\bf 1}_j$ with ${\bf 1}_j$ the $j$th
standard unit vector in $\Q^n$, and with the family is associated the
polytope
$$\P^{\F}\quad :=\quad
\conv \,\{\,{\bf 1}_F\,:\,F\in \F\,\}\quad\subset\quad \Q^n\ .$$
Extending $w$ to $\Q^n$ by $w(x):=\sum_{j=1}^n w(j)x_j$, the
problem reduces to maximizing the linear functional $w$ over
$\P^{\F}$. This leads to the study of {\em facets}
of $\P^{\F}$; when these can be suitably controlled, the problem
is polynomial time solvable via the ellipsoid method \cite{Kha}
for linear programming.

In this article we show that the geometric approach can be
usefully exploited so as to yield a widely applicable general
methodology for convex combinatorial optimization as well. Our
framework leads to the study of {\em edge-directions} of
$\P^{\F}$; when these can be suitably controlled, the problem is
efficiently solvable via zonotope (or hyperplane arrangement)
methods as follows.

\bt{Main}
For any fixed $d$, there is a strongly polynomial oracle time algorithm
solving convex combinatorial optimization over any edge-guaranteed family
presented by a membership oracle.
\et
The precise definition of an edge-guaranteed family will be given
in Section 2: all families underlying the various applications in Section 3
naturally possess this property. The assumption of fixed $d$ is also natural
and necessary: already for $d=1$, the problem generalizes linear
combinatorial optimization which is frequently intractable; and when $d$ is
variable, the problem captures NP-hard instances even for the simple
power set family $\F=2^N$, see Example \ref{Quadratic} below.

The main part of the proof of this theorem is a reduction of the
convex combinatorial optimization problem over a family $\F$ to
the solution of polynomially many standard linear combinatorial
optimization counterparts over the same family $\F$. The reduction
makes use of several results about zonotopes which are available
in combinatorial and computational geometry, but have not been so
far integrated and harnessed in a systematic way to discrete optimization.

The repeated solution of each of the standard linear combinatorial
optimization counterparts can be done following either one of the
following two approaches. The first is to use any efficient ad-hoc
algorithm available in the literature for the specific family $\F$
at hand: this typically leads to the best overall running time.
This approach indeed applies to all of the applications discussed
in Section 3, since each admits a very fast ad-hoc algorithm
(ranging from simple greedy to sophisticated min-cost flow). The
second approach, which is generic and works for any $\F$, takes
advantage of the fact that a {\em test set} (cf. \cite{ST}) is
readily available for any edge-guaranteed family, and (see
\cite{SW} and references therein), using scaling \cite{EK} and
Diophantine approximation \cite{FT}, allows the efficient solution
of the necessary linear optimization counterparts.

\vskip.4cm

The article is organized as follows. In Section 2 we prove Theorem
\ref{Main} as well as some other results, and discuss some
relevant issues, as follows. In \S 2.1 we discuss the necessary
preliminaries on zonotopes and edge-directions. In \S 2.2 we prove
Theorem \ref{Reduction} providing the reduction of convex
combinatorial optimization to polynomially many standard linear
combinatorial optimization counterparts. In \S 2.3 we discuss the
generic approach for solving the counterparts and combine it with
Theorem \ref{Reduction} to establish Theorem \ref{Main}. In \S 2.4
we discuss the problem of finding short monotone paths on
$(0,1)$-polytopes, provide Lemma \ref{Long} which is a certain
$(0,1)$-analog of the Klee-Minty cube \cite{KM}, and raise some
questions. In \S 2.5 we consider classes of edge-guaranteed
families and conclude Corollary \ref{Class} concerning such
{\em edge-well-behaved classes} (defined therein). In \S 2.6 we
discuss projection representation which sometimes helps control
edge-directions. Section 3 is devoted to applications:
in \S 3.1 we discuss quadratic assignment and matroids, whereas in
\S 3.2 we make use of projection representation and discuss in
detail the broadly applicable shape vector partitioning problem.
The last Section 4 contains some final remarks and open problems.

\vskip.4cm

We conclude this introduction with some comments. First, our results make
use of and provide an efficient enumeration of the vertices of the polytope
$\P^{\F}_w:=\conv \{w(F): F \in \F\}$ which is a projection of $\P^{\F}$;
as the maximum of a convex functional $c$ over $\P_w^{\F}$ is attained at a
vertex and each vertex has the form $w(F)$ with $F\in\F$, this provides a
strategy for addressing the convex combinatorial optimization problem.
One of the difficulties we overcome is that the number of sets in $\F$
is typically exponential in $n$ and hence it is generally impossible to
construct $\P^\F_w$ directly efficiently. As a consequence of our
efficient vertex enumeration of $\P^{\F}_w$, our results immediately
extend to the larger class of problems where $c$ is any functional
which is guaranteed to attain a maximum over $\P_w^{\F}$ at a vertex, e.g.,
when $c$ is (edge-)quasi-convex on $\P_w^{\F}$, see \cite{HR}.
In particular, our results extend to
functions $c$ which are (asymmetric) Schur convex when the
edge-directions of $\P_w^{\F}$ are differences of standard unit
vectors in $\Q^d$, again see \cite{HR}.  Also, the results can be
generalized to some extent from $(0,1)$-problems to integer
programming. Second, we note that in studying edge-directions, we
make use of projections of polytopes; thus the Billera-Sturmfels
theory of fiber polytopes \cite{BSt}, related to various aspects
of polytope projection, may be helpful in the classification of
edge-well-behaved classes of families. Also, some new questions
that we raise about graphs of $(0,1)$-polytopes might be addressed
through the Billera-Sarangarajan universal embedding of such
polytopes as travelling salesman polytopes \cite{BSa}.

We hope this exposition will make our framework a widely
accessible tool in the arsenal of discrete optimization, and will
stimulate the study of edge-directions of polytopes $\P^{\F}$ for
various combinatorially defined families. Since convex
combinatorial optimization is often intractable, there is also
much room for the study of approximation algorithms for this
problem for various families $\F$, and we hope this article will
stimulate research on this yet unexplored ground.

\section{Edge-Directions and the Algorithmic Solution}

\subsection{Edge-directions and zonotopes}

We start by introducing the necessary terminology and collecting
several facts that we shall make use of; for some we only provide
a reference and for others we provide a short proof.

The {\em zonotope} generated by a set of vectors $E=\{e^1,\dots,e^m\}$ in
$\Q^d$ is the Minkowsky sum
$$\Z\,=\,\zone(E)\,:=\,\sum_{i=1}^m\ [-e^i,e^i]
\,=\,\left\{\sum_{i=1}^m \lambda_i e^i\,:\,-1\leq\lambda_i\leq
1\right\} \,=\,\conv\left\{\sum_{i=1}^m \lambda_i e^i\,:\,
\lambda_i=\pm 1\right\} \,\subset\,\Q^d \ .$$ The following bound
on the number of vertices of zonotopes has been rediscovered many
times over the years; see e.g. \cite{Buc,Har} for early references
and \cite{GS,Zas} for recent extensions and refinements.
\bl{Zonotope1}
The number of vertices of any $d$-dimensional
zonotope generated by $m$ vectors is at most
$2\sum_{i=0}^{d-1}{{m-1}\choose i}$. Thus, for fixed $d$ it is
$O(m^{d-1})$ and hence polynomially bounded in $m$.
\el
Each vector $a\in\Q^d$ is also interpreted as the linear functional on
$\Q^d$ given via the standard inner product $a\cdot x=\sum_{i=1}^d
a_i x_i$. The {\em normal cone} of a polytope $\P$ at its face $F$
is the (relatively open) cone $\C_\P^F$ of those linear
functionals $a$ which are maximized over $\P$ precisely at points
of $F$. The following computational analogue of Lemma
\ref{Zonotope1} is provided by the algorithm in \cite{EOS,ESS}
(the latter reference provides a necessary correction of the
former); some extensions are again in \cite{GS}.
\bl{Zonotope2}
Fix any $d$. Then all vertices of any $d$-dimensional zonotope
$\Z$ generated by $m$ given vectors can be listed, each vertex $u$
along with a linear functional $a(u)\in\C_\Z^u$ uniquely maximized
over $\Z$ at $u$, in strongly polynomial time using $O(m^{d-1})$
arithmetic operations.
\el
Note that throughout we are mainly interested in {\em strongly polynomial
(oracle) time  algorithms}, that is, algorithms that perform a polynomial
number of arithmetic operations (and calls to the relevant oracles if any)
and are also polynomial time in the Turing computation model.

\vskip.4cm
The collection of normal cones of a polytope $\P$ at all faces is called the
{\em normal fan} of $\P$ (see \cite{Gru}). A polytope $\Z$ is a {\em refinement} of a polytope $\P$ if the closure of each normal cone
of $\P$ is the union of closures of normal cones of $\Z$. A standard result
shows that $\Z$ is a refinement of $\P$ if and only if the normal cone of every
vertex of $\Z$ is contained in the normal cone of some vertex of $\P$, and we
will use this property interchangeably with the above definition of refinement.
A {\em direction} of an edge $[u,v]$ of a polytope $\P$ is any
nonzero scalar multiple of $v-u$. We provide a simple proof of the following
fact (cf. \cite{GS}) which is quite central to our considerations.

\bl{Refinement}
Let $\P$ be a polytope and let $E$ be a finite set of vectors containing
a direction of every edge of $\P$. Then the zonotope $\Z:=\zone(E)$
generated by $E$ is a refinement of $\P$.
\el
\bpr
Let $E=\{z^1,\dots,z^m\}$. Consider any vertex $u$ of $\Z$.
Then $u=\sum_{i=1}^m \lambda_i z_i$ for some $\lambda_i=\pm1$
and hence its normal cone $\C^u_\Z$ consists of those $a$ satisfying
$a\cdot\lambda_i z_i>0$ for all $i$. Let $v$ be a vertex of $\P$ at which
some such $\hat a$ (belonging to $\C^u_\Z$) is maximized over $\P$.
Consider any edge $[v,w]$ of $\P$.
Then $v-w=\alpha_i z_i$ for some scalar $\alpha_i\neq 0$ and some $z_i$,
and $0\leq{\hat a}\cdot(v-w)={\hat a}\cdot\alpha_i z_i$, implying
$\alpha_i\lambda_i>0$. It follows that every $a$ in the cone $\C^u_\Z$ of the
vertex $u$ of $\Z$ satisfies $a\cdot(v-w)>0$ for every edge of $\P$ containing
$v$ and therefore $a$ is also in the cone $\C^v_\P$ of the vertex $v$ of $\P$,
and hence $\C^u_\Z\subseteq\C^v_\P$. Since $u$ was arbitrary, it follows
that the normal cone of every vertex of $\Z$ is contained in the normal cone
of some vertex of $\P$ and we are done by the aforementioned standard result.
\epr

Finally, we need the following statement about
edge-directions of linear images of polytopes.

\bl{Projection}
Let ${\cal Q}:=\omega(\P)$ be the image of a polytope $\P$ under a
linear map $\omega$. Then every direction $q$ of an edge of $Q$ is
the image under $\omega$ of some direction $p$ of an edge of $\P$.
\el
\bpr
Let $q$ be a direction of an edge
$[x,y]$ of ${\cal Q}$. Consider the face $F:=\omega^{-1}([x,y])$ of
$\P$. Let $V$ be the set of vertices of $F$ and let $U=\{u\in
V\,:\,\omega(u)=x\,\}$; as $x \neq y$, $U\neq V$.
Further, as the graph of $F$ is connected there must be an edge $[u,v]$
of $F$, and hence of $\P$, for some $u\in U$ and $v\in V\setminus U$.
Then $\omega(v)\in(x,y]$ hence $\omega(v)=x+\alpha q$ for some
$\alpha\neq 0$. Therefore
$q={1\over\alpha}(\omega(v)-\omega(u))=\omega({1\over\alpha}(v-u))=\omega(p)$
with $p:={1\over\alpha}(v-u)$, a direction of the edge $[u,v]$ of $\P$.
\epr

\subsection{Reduction of convex to linear combinatorial optimization }

We now reduce convex to linear combinatorial optimization. We make
the following assumptions. The ground set is $N:=\{1,\dots,n\}$
and the family $\F\subseteq 2^N$ is {\em edge-guaranteed}, which
means that it is nonempty and comes with an explicit set
$E=\{e^1,\dots,e^m\}\subseteq\Q^n$ of vectors guaranteed to
contain a direction of each edge of the polytope $\P^{\F}=\conv
\{\,{\bf 1}_F\,:\,F\in \F\,\}$ associated with $\F$. In this
subsection we assume that $\F$ is presented by a {\em linear
combinatorial optimization oracle} that, given
$b:N\longrightarrow\Q$, returns a family member $F\in\F$ of
maximum weight $b(F)=\sum_{j\in F}b(j)$. The convex functional
$c:\Q^d\longrightarrow \Q$ is presented by an {\em evaluation oracle}
that, given $x\in\Q^d$, returns the value $c(x)$. The
weighting $w:N\longrightarrow\Q^d$ is given by an explicit list
$w(1),\dots,w(n)\in\Q^d$. We consider $d$ as fixed; otherwise, as
mentioned before, the problem becomes intractable at once even for
the simple family $\F=2^N$, see Example \ref{Quadratic} below. The
following algorithm, applied to the data above, provides a
reduction of convex to linear combinatorial optimization.

\ba{Algorithm}
Given data as above, perform the following steps:
\begin{enumerate}
\item
Consider the linear map $\omega:\Q^n\longrightarrow\Q^d$ defined
by $\omega(x):=\sum_{j=1}^n w(j)x_j\ $.
\begin{enumerate}
\item
Compute the image $\omega(E)=\{\omega(e^1),\dots,\omega(e^m)\}$
of $E$ under $\omega$.
\end{enumerate}
\item
Consider the zonotope
$\Z:=\zone(\omega(E))=\sum_{i=1}^m [-\omega(e^i),\omega(e^i)]$
in $d$-space $\Q^d$\ .
\begin{enumerate}
\item
Compute the list $\{u^1,\dots,u^k\}$ of all vertices of $\Z$.
\item
For each $u^i$ compute a linear functional
$a^i\in\C_\Z^{u^i}$ in the normal cone of $\Z$ at $u^i$.
\end{enumerate}
\item
\begin{enumerate}
\item
For each $a^i$ compute $b^i:N\longrightarrow\Q$
defined by $b^i(j):=a^i\cdot w(j)=\sum_{t=1}^d a^i_t w(j)_t\ $.
\item
For each $b^i$ query the oracle of $\F$ and obtain
$F^i\in\F$ of maximum weight $b^i(F^i)$.
\item
For each $F^i$ query the oracle of $c$ and obtain
the value $c(w(F^i))=c(\sum_{j\in F^i}w(j))$.
\end{enumerate}
\item
Output $F^i\in\F$ of maximum value $c(w(F^i))$ among $F^1,\dots,F^k$.
\end{enumerate}
\ea

\bt{Reduction}
Algorithm \ref{Algorithm} solves the convex
combinatorial optimization problem with data as above in strongly
polynomial oracle time using $O(nm^{d-1})$ arithmetic operations
and $O(m^{d-1})$ queries of the linear combinatorial
optimization oracle of $\F$ and the evaluation oracle of $c$.
\et
\bpr
First we justify the algorithm. Recall the polytope
$\P^\F_w=\conv \{\,w(F)\,:\,F\in \F\,\}$. As
$$\P^\F_w \ =\ \conv\{\,w(F)\,:\,F\in \F\,\}
\ =\ \conv\{\,\omega({\bf 1}_F)\,:\,F\in \F\,\} \ =\
\omega(\P^\F)\ ,$$
$\P^\F_w$ is the image of $\P^\F$ under the
linear map $\omega$ defined in step 1 of the algorithm.  Thus, by
Lemma \ref{Projection}, the image $\omega(E)$ of $E$ under
$\omega$ contains a direction of every edge of $\P^\F_w$.
Therefore, by Lemma \ref{Refinement}, the zonotope
$$\Z \quad := \quad \zone(\omega(E))\quad=\quad
\sum_{i=1}^m [-\omega(e^i),\omega(e^i)]$$ defined in step 2 is a
refinement of $\P^\F_w$. Now consider any vertex $v$ of $\P^\F_w$.
Since $\Z$ refines $\P^\F_w$, the normal cone of $\P^\F_w$ at $v$
contains the normal cone of $\Z$ at some vertex $u^i$ of $\Z$
found in step 2a. This implies that the corresponding linear functional
$a^i\in\C^{u^i}_\Z$ found in step 2b is maximized
uniquely over $\P^\F_w$ at $v$.
Now, consider the corresponding weighting $b^i$ defined in step 3a.
As $v$ is the unique maximizer of $a^i$ over
$\P^\F_w=\conv\{\,w(F)\,:\,F\in \F\,\}$, we have
$$b^i(F)\ =\ \sum_{j\in F} a^i\cdot w(j)
\ = \ a^i\cdot \sum_{j\in F}w(j)\ = \ a^i\cdot w(F)\ \leq\ a^i\cdot v$$
for each $F\in\F$, with equality if and only if $w(F)=v$.
Thus, the member $F^i\in\F$ obtained in step
3b from the linear combinatorial optimization oracle of $\F$ when
maximizing $b^i$ has $v=w(F^i)$. It follows that every vertex of
$\P^\F_w$ equals $w(F^i)$ for some $F^i$ obtained in step 3b.
Since $c$ is convex, the maximum value $c(w(F))$ of $F\in\F$
occurs at some vertex $v=w(F^i)$ of $\P^\F_w$. Thus the member
$F^i\in\F$ output by the algorithm in step 4, which has maximum
value $c(w(F^i))$ among the values computed in step 3c, is an
optimal solution to the convex combinatorial optimization problem.

Next we verify the claimed complexity, where, as explained, $d$ is
considered fixed. The computation of the linear image $\omega(E)$
in step 1a takes $O(dnm)=O(nm)$ operations. By Lemma
\ref{Zonotope1}, the number of vertices of the zonotope $\Z$
defined in step 2 satisfies $k=O(m^{d-1})$, and the computation of
these vertices $u^i$ and of corresponding linear functionals $a^i$
in steps 2a and 2b requires $O(m^{d-1})$ operations by Lemma
\ref{Zonotope2}. The number of queries in step 3b of the oracle of
$\F$ and in step 3c of the oracle of $c$ are $k=O(m^{d-1})$ as
claimed. The computation of each $b^i$ in step 3a and of each
$w(F^i)$ in step 3c take $O(dn)=O(n)$ operations totalling
together over all $i$ to $O(kn)=O(nm^{d-1})$ arithmetic
operations. Finally, the arithmetic complexity of finding the
maximum among the $k$ values $c(w(F^i))$ in step 4 is
$O(k)=O(m^{d-1})$. Thus, the dominant arithmetic complexity is
$O(nm^{d-1})$ as claimed. \epr

As the proof shows, convex combinatorial optimization is solved by
enumerating all vertices of the polytope $\P^\F_w$ and picking the
best. While each vertex of $\P^\F_w$ is the image under $w$ of
some vector {\bf 1}$_F$, with $F \in \F$, the difficulty is that
the number of sets in $\F$ is typically exponential in $n$ and
hence it is generally impossible to construct $\P^\F_w$ directly
in polynomial arithmetic complexity (in particular, each {\bf 1}$_F$
is a vertex of $\P^\F$). The efficient construction of
$\P^\F_w$ is made possible by the given set of edge-directions of
$\P^\F$ and by proceeding, indirectly, through the zonotope $\Z$
that refines $\P^\F_w$. While the number of vertices of $\Z$ can
be much larger than that of $\P^\F_w$, the vertices of $\Z$ can be
better controlled and this leads to the polynomial complexity
bound. So if the polytope $\P^\F$ of a family $\F$ admits a
relatively small set containing a direction of each edge which can
be efficiently constructed or even characterized, then the problem
is efficiently reducible.

\subsection{Generic solution of the linear
combinatorial optimization counterparts }

We now discuss how to realize an oracle that will repeatedly solve
each of the linear combinatorial optimization counterparts queried
upon in Algorithm \ref{Algorithm}. As before, we assume that our
family $\F$ is edge-guaranteed and hence nonempty and comes with
an explicit set $E=\{e^1,\dots,e^m\}$ containing a direction of
each edge of $\P^\F$. We assume moreover that the family comes
with one member $F^0\in\F$ to start with. We consider the
following three oracle presentations of $\F$.
\begin{itemize}
\item
{\em Membership oracle}: when queried about $F\subseteq N$, this oracle
asserts whether or not $F\in \F$.
\item
{\em Augmentation oracle}: when queried about $F\in\F$ and
$b:N\longrightarrow\Q$, this oracle returns a family member ${\hat F}\in\F$
with $b({\hat F})>b(F)$ or asserts that $F$ has maximum weight in $\F$.
\item
{\em Linear combinatorial optimization oracle}: when queried about
$b:N\longrightarrow\Q$, this oracle returns a family member
$F\in\F$ of maximum weight $b(F)$.
\end{itemize}

\bl{Membership}
For any edge-guaranteed family, a membership oracle enables to simulate
an augmentation oracle in strongly polynomial oracle time.
\el
\bpr
Without loss of generality,
assume that each $e^i$ is a $\{-1,0,1\}$-vector.  The simulation
is simple. Consider a query about $F$ and $b:N\longrightarrow\Q$.
Call an edge-direction $e^i$ {\em improving} if
$b(e^i)=\sum_{j=1}^n b(j)e^i_j>0$; call it {\em admissible at} $F$
if ${\bf 1_F}+e^i$ is a $\{0,1\}$-vector whose support
$F^i:=\mbox{supp}({\bf 1_F}+e^i)$ is in $\F$. If there is an
edge-direction $e^i$ which is both improving and admissible then
return $F^i$; otherwise assert that $F$ has maximum weight in
$\F$. The simulation works correctly since, as is well known, a
vertex $u$ is not a maximizer of a linear functional $b$ over a
polytope if and only if the polytope has an edge $[u,v]$ for some
vertex $v$ with $b(v)>b(u)$. \epr The next lemma is from
\cite{GL,SWZ}; see \cite{SW} for the state of the art on this line
of research. It involves a computationally heavy Diophantine
approximation step \cite{FT} and a scaling step \cite{EK}. We
include an outline of the proof, which is relevant for the
discussion in the next subsection.
\bl{Augmentation}
For any family, an augmentation oracle enables to simulate a linear
combinatorial optimization oracle in strongly polynomial oracle time.
\el
\bpr
We outline the simulation. Consider query about
$b:N\longrightarrow\Q$. The Diophantine approximation step (see
\cite[Theorem 3.3]{FT}) replaces $b$ by
$a:N\longrightarrow\mathbb{Z}$ with the following two properties:
first, it is equivalent to $b$ in that, for any pair $F,G\subseteq
N$, it satisfies $a(F)\leq a(G)$ if and only $b(F)\leq b(G)$; and
second, the maximum number of bits $k:=1+\max_{j\in N}\lceil\log
|a(j)|\rceil$ in the binary representation of the weight under $a$
of any element is polynomial in $n$.

The scaling step (inspired by \cite{EK}) is the following. Applying
a simple transformation (see \cite{SWZ}) we may assume $a$ is nonnegative.
Following the prof of \cite[Theorem 9.2]{GL},
for $i=0,\dots,k$, starting with $F^i\in\F$, find a maximizing $F^{i+1}\in\F$
with respect to the weighting $a^i:=\lfloor 2^{i-k}a\rfloor$. As shown
in \cite{GL} or \cite{SWZ}, each $F^{i+1}$ is obtained from $F^i$ by calling
the augmentation oracle at most $n$ times. Thus, the maximizer $F^{k+1}$ of
$a=a^k$ is found using at most $kn$ calls, and since $k$ is polynomial in $n$,
the desired maximizer of $b$ is obtained in strongly polynomial oracle time.
\epr
We can now prove Theorem \ref{Main}. As discussed before, the family
$\F\subseteq 2^N$ comes with one explicit $F^0\in\F$ and an explicit
set $E=\{e^1,\dots,e^m\}\subset\Q^n$ containing a direction of each
edge of $\P^\F$, and is presented by a membership oracle.
The complexity is measured in terms of $n$ and $m$.

\vskip.4cm\noindent{\bf Theorem \ref{Main}}
{\em For any fixed $d$, there is a strongly polynomial oracle time
algorithm solving convex combinatorial optimization over any
edge-guaranteed family presented by a membership oracle.}

\vskip.4cm
\bpr
Theorem \ref{Reduction} guarantees that Algorithm
\ref{Algorithm} solves the convex combinatorial optimization
problem over $\F$ efficiently using a linear combinatorial
optimization oracle which, by Lemmas \ref{Membership} and
\ref{Augmentation}, can be efficiently simulated from the
membership oracle presenting $\F$.
\epr
We conclude this subsection with several important remarks which lead
to the discussion in the next subsection. First, the complexity behind
Theorem \ref{Main} is quite horrendous: for each linear
combinatorial optimization counterpart invoked by Algorithm
\ref{Algorithm}, an application of the Diophantine approximation
step which takes $O(n^8)$ arithmetic operations \cite{FT} is required.
However, improved complexity bounds follow from Theorem \ref{Reduction}
when Algorithm \ref{Algorithm} is used with a more efficient linear
combinatorial optimization oracle whenever a particular family admits one.

Second, what about real data and real arithmetic computation (where pairs of
real numbers can be added, multiplied or compared in unit time)?
Algorithm \ref{Algorithm} remains valid and polynomial and the analog of
Theorem \ref{Reduction} (with ``strongly polynomial oracle time" replaced by
``polynomially many real arithmetic operations and queries") holds.
So does the conversion of the membership oracle to the augmentation oracle
manifested by Lemma \ref{Membership} above. However, the proof of
both parts of Lemma \ref{Augmentation} (scaling and Diophantine approximation)
breaks down for real data, and the conversion of the augmentation oracle
to the optimization oracle is no longer available.

Can these obstacles be waved and does the real analog of Theorem \ref{Main}
remain valid? Our families are edge-guaranteed, which is stronger than
having a test set (see \cite{ST} and references therein) and even more so
than having a mere augmentation oracle:
can we take advantage of that and simulate standard linear combinatorial
optimization directly and more efficiently, avoiding scaling and
Diophantine approximation? We discuss some of these issues next.

\subsection{On the Hirsch conjecture and the
Klee-Minty problem for $(0,1)$-polytopes}

A form of the {\em Hirsch conjecture}, open to date, asks whether
the diameter of (the graph of) every $n$-polytope with $f$ facets
is bounded above by a polynomial in $n$ and $f$; for all one
knows, the linear upper bound $f-n$ may suffice, see \cite{Kal,KK}.
Also open is the analogous form of the {\em monotone Hirsch conjecture}
asking whether the shortest increasing path under any linear functional
from any vertex to some maximizing vertex is polynomially
bounded in $n$ and $f$. Both variants are true for $(0,1)$-polytopes
\cite{Nad} as well as for some more general classes of integer
polytopes \cite{DO,KO}. The following slightly stronger form,
relevant to the discussion below, holds.

\bl{Diameter}
Any $(0,1)$ $n$-polytope $P$ admits, under any linear functional, a
nondecreasing path from any vertex $u$ to any maximizing vertex
$w$, of length at most $n$ using no edge-direction twice.
\el
\bpr
The claim being trivial for $n=0$, we proceed by induction. If $u$
and $w$ lie on a common proper face of $P$ then induction takes
over. Otherwise, there is a nondecreasing arc $(u,v)$; pick any
$i$ with $u_i\neq v_i$; then $v_i=w_i$ and hence, by induction,
there is a nondecreasing $(v,w)$-path of length at most $n-1$
using no edge-direction twice on the face $F:=\{x\in
P\,:\,x_i=w_i\}$. As no edge of $F$ can have direction $v-u$, this
path preceded by the arc $(u,v)$ gives the desired $(u,w)$-path.
\epr The {\em effective Hirsch conjecture} asks, broadly, whether
a monotone path could be efficiently traced. To make things
precise, the presentation of the polytope has to be specified. For
instance, tracing such a path in strongly polynomial time for
polytopes presented by linear inequalities would imply a strongly
polynomial time algorithm for linear programming via the simplex
method which does not seem likely; but tracing it in
subexponential time is possible \cite{Kal}. A natural question is:
how long can an {\em arbitrary} increasing path be? A classical
construction by Klee and Minty \cite{KM} transforms the $n$-cube
so as to admit increasing paths of exponential length $2^n$. But
$(0,1)$-polytopes are very special, as shows Lemma \ref{Diameter}.
How long, then, can an arbitrary increasing path in a
$(0,1)$-polytope be? We now show that, unfortunately,
$(0,1)$-polytopes admit such paths of length exponential in
the dimension as well: in this sense, the following lemma can be
regarded as a $(0,1)$ analog of the Klee-Minty cube; we thank Tal
Raviv for a related discussion.
\bl{Long}
For every $n$ there is a $(0,1)$-polytope of dimension less than
${1\over 4}n^4$ with $n!$ vertices that admits a Hamiltonian
(and hence $n!$-long) nondecreasing path under every linear functional.
\el
\bpr
The {\em Young polytope} $Y_{n-2,2}$ is the convex hull of all
${n\choose 2}\times{n\choose 2}$ matrices of permutations of edges
of the complete graph $K_n$ induced by the $n!$ permutations of
its vertices. For instance, the matrix corresponding to the
permutation of vertices $\sigma=(1, 2, 3, 4)$ (in cycle notation)
is
$$\small{
\Sigma\quad =\quad \bordermatrix{& 12 & 13 & 14 & 23 & 24 & 34 \cr
12 & 0 & 0 & 0 & 1 & 0 & 0 \cr 13 & 0 & 0 & 0 & 0 & 1 & 0 \cr 14 &
1 & 0 & 0 & 0 & 0 & 0 \cr 23 & 0 & 0 & 0 & 0 & 0 & 1 \cr 24 & 0 &
1 & 0 & 0 & 0 & 0 \cr 34 & 0 & 0 & 1 & 0 & 0 & 0 \cr } \quad.}$$
The claims about the dimension and number of vertices of
$Y_{n-2,2}$ are obvious. In \cite{Onn1} it was shown that
$Y_{n-2,2}$ is $2$-neighborly, that is, its graph is the complete
graph $K_{n!}$: the very existence of $2$-neighborly
$(0,1)$-polytopes is an amazing fact in itself! It follows that if
$a$ is any linear functional, then any ordering $v_1,\dots,v_{n!}$
of the vertices of $Y_{n-2,2}$ satisfying $a(v_1)\leq\cdots\leq
a(v_{n!})$ gives a Hamiltonian (and hence exponentially long)
nondecreasing path under $a$. \epr Here, however, we are
especially interested in the polytopes $\P^\F$ of {\em
edge-guaranteed} families. Lemmas \ref{Membership} and
\ref{Augmentation} imply that for such a family, a monotone path
can be traced in time strongly polynomial in $m$ and $n$, alas,
for {\em rational} functionals only, and using the heavy
Diophantine approximation procedure. Can we do better? what is the
maximal length $I(n,m)$ of any increasing path in any
$n$-dimensional $(0,1)$-polytope with $m$ pairwise nonproportional
edge-directions?

\subsection{Edge-well-behaved classes}

In most applications, in particular all of those discussed in
Section 3, one is concerned with a class of families possessing
some unifying structure. It is therefore useful to make some
formal definitions regarding such classes and then use it to
obtain a suitable corollary of Theorem \ref{Main}. For $n\geq 0$
let $N=\{1,\dots,n\}$ as before and let ${\cal U}_n$ be the set of
all families with ground set $N$,
$${\cal U}_n \quad:=\quad 2^{2^{\{1,\dots,n\}}}
\quad=\quad\{\F\ :\ \F\subseteq 2^N\}\ .$$
A {\em class of families} is a (typically infinite) set of families
$\C=\biguplus_{n\geq 0}\C_n$ with $\C_n\subseteq {\cal U}_n$ for
all $n$.
\bd{Definition}
A class $\C$ is {\em edge-well-behaved} if there is a
polynomial time algorithm that, given $n$, produces a set
$E_n=\{e^1,\dots,e^{m(n)}\}\subseteq\{-1,0,1\}^n$ with
respect to which every $\F\in \C_n$ is edge-guaranteed. In
particular, $m(n)$ is polynomial in $n$ and each $\F\in\C$ is nonempty.
\ed
While the existence of such a ``uniform" polynomial
time algorithm that produces sets containing edge-directions for
the polytopes of all families in the class may seem a strong
assumption, we will see in Section 3 that in many applications
such an algorithm is readily available. Also, the assumption that the
edge-directions are $(-1,0,1)$-valued is not restrictive since, \break
for $(0,1)$-polytopes, each edge is a difference of
two vertices hence admits a $(-1,0,1)$-direction. The next
corollary follows at once from Theorem \ref{Main}; here the
complexity is in terms of $n$ only.
\bc{Class}
Fix any $d$. Then for any edge-well-behaved class there is a strongly
polynomial oracle time algorithm that solves the convex combinatorial
optimization problem over any family in the class which is
presented by a membership oracle.
\ec
While this statement may seem a reformulation of Theorem \ref{Main},
it is natural and useful in uniformly establishing the polynomial
solvability in all of the applications discussed in Section 3.

\subsection{Projection representation and circuits}

We conclude Section 2 by discussing a useful setup that helps in
controlling edge-directions, and which is used and demonstrated in
the application given in \S 3.2 in the sequel.
 A {\em circuit} of an $r\times s$ matrix $A$ is a nonzero
solution $z\in\Q^s$ of the system $Az=0$ whose support is
inclusion-minimal. It is known (cf. \cite[Exercise 10.14]{Roc})
that any nonzero solution of $Az=0$ has a {\em
conformal circuit decomposition}, i.e. can be expressed as
$z=\sum_i \alpha_i z^i$ with the $\alpha_i$ positive scalars and
the $z^i$ pairwise nonproportional circuits such that $ z^i_j z_j
>0$ for all $i$ and all $j\in\mbox{supp}(z^i)$. Consider the
{\em standard polytope} $\P=\{x\in\Q^s\,:\,Ax=b,\ l\leq x\leq u\}$
defined by $A$, right-hand side $b\in\Q^r$, lower bound
$l\in\Q^s$, and upper bound $u\in\Q^s$. We provide a short proof
of the following useful property of edge-directions of the
standard polytope (see \cite{ORT03} for a refinement of this
property which characterizes edge-directions).

\bl{Circuits}
Each edge-direction of a standard $\P$ is a circuit of its defining matrix $A$.
\el
\bpr
Consider any $x,y\in\P$. Then $A(y-x)=0$ so $y-x$ admits a
conformal circuit decomposition $y-x=\sum_i \alpha_i z^i$. It is
then not hard to verify that for every circuit $z^i$ participating
in that decomposition, both $x+\alpha_i z^i$ and $y-\alpha_i z^i$
satisfy the lower and upper bounds and hence are in $\P$. They
belong, moreover, to any face $F$ containing both $x$ and $y$.
Indeed, pick any $c$ in the normal cone $\C_\P^F$: then $c\cdot
x\geq c\cdot(x+\alpha_i z^i)$ which implies $c\cdot z^i \leq 0$,
and $c\cdot y\geq c\cdot(y-\alpha_i z^i)$ which implies $c\cdot
z^i \geq 0$. It follows that $c\cdot z^i =0$ and hence $c\cdot
(x+\alpha_i z^i)=c\cdot x$ and $c\cdot (y-\alpha_i z^i)=c\cdot y$
implying $x+\alpha_i z^i, y-\alpha_i z^i\in F$. Now, if the
decomposition $y-x=\sum_i \alpha_i z^i$ involves more than one
circuit, say $z^1,z^2$, then any face containing $x,y$ contains
the three non-collinear points $x,x+\alpha_1 z^1,x+\alpha_2 z^2$
and hence is not an edge. So if $[x,y]$ is an edge of $\P$ then
$y-x=\alpha z$ for some circuit $z$. Any direction of that edge is
a nonzero multiple of $z$ and hence a circuit of $A$. \epr Lemma
\ref{Circuits} implies that any inclusion-maximal set $Z$ of
pairwise nonproportional circuits of $A$ contains a direction of
each edge of $\P$. So if the size of $A$ is $r\times s$ then $\P$
admits such a set $Z$ with no more than ${s\choose r}$ elements.
Every polytope ${\cal Q}$ is the linear image of a standard polytope:
if ${\cal Q}=\{x\,:\, Bx\leq b \}$ is a description by inequalities then,
adding a suitable ``slack" vector $y$, we get ${\cal Q}=\phi(\P)$ with
$\P=\{(x,y)\,:\,Bx+Iy=b,\ y\geq 0\}$ and with $\phi$ the ``$y$
forgetting" projection $\phi(x,y)=x$. In particular, the polytope
$\P^\F$ of any family $\F$ is the linear image $\P^\F=\phi(\P)$ of
a standard polytope $\P=\{x\in\Q^s\,:\,Ax=b,\ l\leq x\leq u\}$.
Typically the number of inequalities describing $\P^\F$ is
exponentially large and hence so is the dimension of $\P$, but
when $\P$ has small dimension, we can benefit from such a
``projection representation" in two ways as follows.

First, if the defining matrix $A$ admits an efficiently
determinable set $Z=\{z^1,\dots,z^m\}\subset\Q^s$ containing a
scalar multiple of each circuit of $A$ then, by Lemmas
\ref{Projection} and \ref{Circuits}, its image $E:=\phi(Z)$
contains a direction of each edge of $\P^\F$, making $\F$ an
edge-guaranteed family. Second, linear combinatorial optimization
over $\F$ can be ``lifted" to linear programming over the polytope
$\P$, giving a way alternative to \S 2.3 for solving the
counterparts called upon by Algorithm \ref{Algorithm}.

\section{Some Applications}

\subsection{Some direct applications}

Here we give two examples where the set of edge-directions can be
directly determined and used.

\be{Quadratic}{\bf Positive semidefinite quadratic assignment.}
The quadratic assignment problem is the following: given a real $n\times n$
matrix $M$, find $x\in\{0,1\}^n$ maximizing the quadratic form
$x^T M x$ induced by $M$; see \cite{PRW} for an overview
of this problem and its applications.
We consider the instance where $M$ is positive semidefinite, in which case
it can be assumed to be presented as $M=W^T W$ with $W$ a given $d\times n$
matrix. If the rank $d$ of $W$ and $M$ is variable then this problem is
NP-hard \cite{HHPR}. For fixed $d$ it is polynomial time solvable \cite{AFLS}.

When $W$ is rational, the problem can be modelled as convex
combinatorial optimization with the following data: the family is
the entire power set $\F=2^N$ of $N$ with the natural
correspondence $\F\leftrightarrow\{0,1\}^n$; the weight of $j\in
N$ is the $j$th column $w(j):=W^j$ of the matrix $W$; and
$c:\Q^d\longrightarrow\Q:x\mapsto ||x||^2$ is the squared standard $l_2$ norm.
Indeed, for each $F\in\F$ we then have ${\bf 1}_F^T W^T W{\bf 1}_F=c(w(F))$.

Now, the polytope of the family $\F$ here is just the $n$-cube
$\P^\F=[0,1]^n$; therefore the trivially computable set of $n$
standard unit vectors $E:=\{{\bf 1}_1,\dots,{\bf 1}_n\}$ contains
a direction of each edge. Thus, the class of all such families is
edge-well-behaved with $m(n)=n$ and Corollary \ref{Class} applies
and guarantees the efficient solution. Here, one obtains a faster
solution by using Algorithm \ref{Algorithm} together with a linear
combinatorial optimization oracle realized by simple sign checking
as follows: given $b:N\longrightarrow\Q$, a member $F\in\F$
maximizing $b(F)$ is simply $F:=\{j\,:\,b(j)>0\}$. \ee

\be{Matroid}{\bf Convex matroid optimization.}
This problem is the special case of convex combinatorial optimization where
$\F$ is either the collection $\cal B$ of bases (considered in \cite{Onn2})
or the collection $\cal I$ of independent sets of a matroid over $N$.
It generalizes classical matroid optimization, first studied in
\cite{Edm2}, and has a rich modelling power on its own: useful matroids
include the forest matroid of a graph and, more generally, the matroid of
linear dependencies of a matrix over a field. For us here it suffices
that the matroid is presented by a membership oracle for $\F$.

It can be derived from the matroid-bases-axioms that the trivially
computable set $D:=\{{\bf 1}_i-{\bf 1}_j\,:\,1\leq i < j \leq n\}$
of ${n\choose 2}$ differences of unit vectors contains a direction
of each edge of the polytope $\P^{\cal B}$ of the family $\cal B$
of bases. Likewise, it can be derived from the
matroid-independence-axioms that the ${n+1}\choose2$-element
union $D\cup E$ of $D$ and the set $E:=\{{\bf 1}_1,\dots,{\bf 1}_n\}$
of unit vectors contains a direction of each edge of the
polytope $\P^{\cal I}$ of the family $\cal I$ of independent sets.
Thus, the class of all such families ${\cal B}$ (respectively,
families ${\cal I}$) is edge-well-behaved with $m(n)={n\choose 2}$
(respectively, $m(n)={{n+1}\choose2}$), and Corollary \ref{Class}
applies and guarantees the efficient solution.

Here too, one obtains a faster solution by using Algorithm
\ref{Algorithm} together with a linear combinatorial optimization
oracle over $\F={\cal B}$ or $\F={\cal I}$ which is efficiently
realizable from a membership oracle for $\F$ using the classical
greedy algorithm (cf. \cite{CCPS,Edm2}) that, given
$b:N\longrightarrow \Q$ makes use of sorting the values $b(j)$ to
find the lexicographically $b$-largest member $F\in\F$ which can
be shown to be the one maximizing $b(F)$. \ee

\subsection{Shaped vector partitioning}

The {\em shaped partition problem} concerns the partitioning of a multiset
$V=\{v^1,\dots,v^n\}$ of $n$ vectors in $d$-space into $p$ parts
so as to maximize an objective function which is convex on the
sum of vectors in each part, subject to constraints on the number
of elements in each part. To describe the problem precisely
we need some notation. A {\em $p$-partition} of the index set
$\{1,\dots,n\}$ of $V$ is an ordered tuple $\pi=(\pi_1,\dots,\pi_p)$
of pairwise disjoint sets whose union is $\{1,\dots,n\}$.
The {\em shape} of a partition is the tuple of cardinalities of its parts,
$|\pi|:=(|\pi_1|,\dots,|\pi_p|)$. In addition to the set of vectors $V$,
the data includes vectors $l,u\in\{0,1\dots,n\}^p$ with $l\leq u$
providing lower and upper bounds on the shape of admissible partitions.
With each partition $\pi$ is associated a $d\times p$ matrix
$$V^{\pi}\quad:=\quad \left[\left(\sum_{i\in\pi_1} v^i\right),
\dots, \left(\sum_{i\in\pi_p} v^i\right)\right]
\quad=\quad \sum_{j=1}^p\left(\sum_{i\in\pi_j} v^i\right){\bf 1}_j^T $$
whose $j$th column is the sum (representing the ``total value") of vectors
assigned to the $j$th part. The data also includes a convex functional
$c:\Q^{d\times p}\longrightarrow\Q$ on $d\times p$ matrices
which ``weighs together" the sums of vectors in the various parts.
The problem is to find a $p$-partition $\pi$ whose shape satisfies the lower and
upper bounds $l\leq|\pi|\leq u$ and which maximizes the value $c(V^{\pi})$.

Shaped partition problems have applications in diverse fields such
as clustering, inventory, reliability, and more - see
\cite{BHR, BH, COR, HOR, HR03, PRW} and references therein.
Here is a typical example.

\be{Clustering}{\bf Minimal variance clustering.}
This is the following problem, which has numerous applications in
the analysis of statistical data: given $n$ sample points
$v^1,\dots,v^n$ in $d$-space, group the points into $p$ clusters
$\pi_1,\dots,\pi_p$ so as to  minimize the sum of cluster variances
$$\sum_{j=1}^p{1\over|\pi_j|}\sum_{i\in \pi_j}
||v^i-({1\over|\pi_j|}\sum_{i \in \pi_j} v^i)||^2\ .$$ We consider
the instance where there are $n=p\cdot m$ points and the
clustering  sought is {\em balanced}, that is, the clusters should
have equal size $m$. Suitable manipulation of the sum of variances
shows that the problem is equivalent to a shaped partition problem
with the lower and upper bounds $l:=u:=(m,\dots,m)$ (forcing the
single shape $|\pi|=(m,\dots,m)$ on partitions), and with the
convex functional (to be maximized) simply as the square of the $l_2$
norm on $d\times p$ matrices, given by $c:\Q^{d\times
p}\longrightarrow\Q:X\mapsto ||X||^2
=\sum_{i=1}^d\sum_{j=1}^p|X_{i,j}|^2$.
\ee

If either the dimension $d$ or the number of parts $p$ is variable, the shaped
partition problem instantly captures NP-hard problems hence is presumably
intractable \cite{HOR}. Therefore, it is interesting to study the worst
case arithmetic complexity in terms of the number $n$ of points with both
$d,p$ fixed. In the special case where there are no shape restrictions
(partitions of all shapes are admissible), an upper bound of $O(n^{d(p-1)-1})$
on the complexity is given in \cite{OSc} and a quite compatible
lower bound of $\Omega(n^{\lfloor {d-1\over2}\rfloor p})$ is in \cite{AvO}.
In the more general case where arbitrary sets of shapes are allowed,
the best upper bound to date is $O(n^{dp^2})$ from \cite{HOR};
while a matching lower bound is unknown, the lower bound $O(n^{d{p\choose 2}})$
from \cite{AlO} on the related number of {\em separable partitions} indicates
that the quadratic term $p^2$ in the exponent may be unavoidable.

\vskip.2cm
We now show how to solve the shaped partition problem efficiently
using our framework. We begin by modelling it as a convex combinatorial
optimization problem. The ground set is taken to be
$N:=\{(i,j)\,:\,1\leq i \leq n,\ 1\leq j \leq p\}$.
Each $p$-partition $\pi=(\pi_1,\dots,\pi_p)$ is encoded as the set
$F_\pi:=\{(i,j)\,:\,i\in\pi_j\}\subseteq N$. The family consists
of all such sets corresponding to $p$-partitions of admissible shapes,
$\F:=\{F_\pi\,:\,l\leq |\pi|\leq u\}$. The weight of element
$(i,j)\in N$ is the $d\times p$ matrix $w(i,j):=v^i {\bf 1}_j^T$
whose $j$th column is $v^i$ and whose other columns are zero.
Finally, the convex functional is simply the given one $c$ defined on
$d\times p$ matrix space. It is not hard to verify that this indeed
casts the shaped partition problem as a
convex combinatorial optimization problem with a ground set
of size $|N|=np$ and weight vectors (matrices) in dimension $dp$.

To show that the class of all such families is edge-well-behaved
we discuss the family polytope $\P^\F\subset\Q^{n\times p}$. The
indicator of a family member $F_\pi\in\F$ is the $(0,1)$-valued
$n\times p$ matrix ${\bf 1}_\pi$ whose $(i,j)$th entry equals $1$
precisely when $i\in\pi_j$. The polytope $\P^\F$ admits a simple
projection-representation as follows. Consider $(n+1)\times p$
matrices whose rows are indexed by $\{0,1,\dots,n\}$. Define lower
and upper bound matrices $L,U$ in terms of the given vectors $l,u$
as follows: for $j=1,\dots,p$ set $L_{i,j}:=0,\ U_{i,j}:=1$ if
$1\leq i\leq n$ and $L_{0,j}:=n-u_j,\ U_{0,j}:=n-l_j$. Let $\P$ be
the {\em transportation polytope} defined by this data, which is
the following standard polytope
$$\P\ :=\ \left\{X\in\Q^{(n+1)\times p}\ :\ \sum_{j=1}^p X_{i,j}=1
\ (1\leq i\leq n), \ \sum_{i=0}^n X_{i,j} =
n\ (1\leq j\leq p),\ L\leq X\leq U\right\}\quad. $$
Then $\P^\F=\phi(\P)$ with
$\phi:\Q^{(n+1)\times p}\longrightarrow\Q^{n\times p}$
the projection erasing the $0$th row of a matrix.

Let $K_{n+1,p}$ be the complete bipartite graph with edge set
$\{\,(i,j)\,:\,0\leq i \leq n,\ 1\leq j \leq p\,\}$ \break
corresponding to this transportation system. Each circuit of
$K_{n+1,p}$ gives an $(n+1)\times p$ matrix supported on that
circuit, with values $\pm 1$ alternating along the edges of the
circuit and $0$ elsewhere. It is well known that each circuit of
the $(n+p)\times (n+1)p$ matrix of coefficients of the equation
system defining $\P$ is proportional to some such
circuit-supported matrix. Let $Z:=\{z^1,\dots,z^m\}$ be the set of
all such $(n+1)\times p$ matrices corresponding to the
$m:=\sum_{i=2}^p {1\over2}{p\choose i}{n+1\choose i}i!(i-1)!$
distinct circuits of $K_{n+1,p}$. Then the projection
$E=\{e^1,\dots,e^m\}:=\phi(Z)$ is a set of $(-1,0,1)$-valued
$n\times p$ matrices which, as explained in \S 2.6, contains a
direction of each edge of $\P^\F$.

Since $p$ is assumed to be fixed, the set of circuits $Z$ and its
projection $E$ are computable in time polynomial in $n$, and
therefore the class of all such families is indeed edge-well-behaved with
$m(n)=\sum_{i=2}^p {1\over2}{p\choose i}{n+1\choose i}i!(i-1)!=O(n^p)$.
Thus, Corollary \ref{Class}
applies and guarantees the efficient solution. As in the examples
in \S 3.1, here too one can obtain a typically faster solution by
using Algorithm \ref{Algorithm} together with a linear
combinatorial optimization oracle realized by lifting to $\P$ and
solving the corresponding transportation problem using available
fast algorithms for bipartite network flows, see \cite{AOST}.
These algorithms, however, are not strongly polynomial and do
depend on the bit size of the bounds $l,u$ and the costs $b^i$
called repeatedly by Algorithm \ref{Algorithm}.

For {\em unrestricted} partitioning, that is, shaped partitioning with
$l=(0,\ldots,0)$ and $u=(n,\ldots,n)$, the characterization of
circuits obtained in \cite{ORT03} shows that circuits of $K_{n+1,p}$
which yield circuits of the matrix of coefficients of the equation
system defining $\P$ correspond to switching a single item from one part
to another. As the number of such circuits is $n{p\choose 2}$,
the class of such families is edge-well-behaved with (improved)
$m(n)=n{p\choose 2}$. Thus, while Theorem \ref{Reduction} with $m(n)=O(n^p)$
and dimension $dp$ implies a complexity bound of $n^{O(dp^2)}$
on the general shaped partition problem, in line with \cite{HOR},
with $m(n)=n{p\choose 2}$ and same dimension $dp$ it implies the
improved bound of $n^{O(dp)}$ on the complexity of the
unrestricted partition problem, in line with \cite{OSc}.

\section{Concluding Remarks}

In this article we have defined the convex combinatorial
optimization problem and shown that it can be solved in strongly
polynomial time for edge-guaranteed families and for
edge-well-behaved classes of families. We have demonstrated
several natural and broad applications that indeed give rise to
edge-well-behaved classes and therefore are efficiently solvable
through our framework.

The polynomial time solvability of linear combinatorial optimization for
{\em facet-well-behaved} classes via the Ellipsoid method \cite{GLS,Kha}
has stimulated over the years a broad body of work on
the identification and characterization of such classes.
A major research program called upon by this paper is an analogous
identification and characterization of {\em edge-well-behaved}
classes of combinatorial families, for which our framework
automatically yields strongly polynomial time solvability of
convex combinatorial optimization.

\vskip.2cm
Some more specific questions are discussed within the
body of our paper, in particular, those in \S 2.3 and \S 2.4
concerning a more efficient generic solution of standard linear
combinatorial optimization over edge-guaranteed families and the effective
Hirsch conjecture for polytopes of such families. What is the
maximal length $I(n,m)$ of any increasing path in any
$n$-dimensional $(0,1)$-polytope with $m$ pairwise nonproportional
edge-directions? Can we trace such a path efficiently while
avoiding scaling and the heavy Diophantine approximation
procedure? Can we trace such a path in polynomially many real
arithmetic operations for real linear functionals?

For solving the standard linear counterparts of a convex combinatorial
optimization problem over a family $\F$ with weighting $w$, one
approach may be to try and augment along edges of the polytope
$\P^\F_w$ downstairs (see discussion following the proof of
Theorem \ref{Reduction}). While a set of edge-directions of that
polytope is available as the projection $\omega(E)$ of the set of
edge-directions of $\P^\F$, this information is not enough: one
needs to know, at any vertex $v$ of $\P^\F_w$, which
edge-directions $\omega(e^i)$ are admissible at $v$, and moreover
- ``how much to walk" - namely, what is the nonnegative scalar
$\alpha$ such that $v+\alpha e^i$ is the new vertex to move to.
When this information is available, even in an abstract setting of
a suitably defined neighborhood oracle, it may be possible to
apply the vertex enumeration methodology of \cite{AF}; and in some
applications, such as unrestricted vector partitioning, this can
be carried out particularly efficiently as in \cite{FOR}.

\vskip.2cm
Finally, as mentioned in the introduction, there is much room for the study of
approximation algorithms for the often intractable convex combinatorial
optimization problem for various classes of families, and we hope this
article will stimulate research on this yet unexplored ground.

\section*{Acknowledgement}

The first author is indebted to Lou Billera for
introducing him to the theory of convex polytopes.

\vskip.5cm\noindent {\small Shmuel Onn}\newline
\emph{Technion - Israel Institute of Technology, 32000 Haifa, Israel}

\emph{email: onn{\small @}ie.technion.ac.il}

\emph{http://ie.technion.ac.il/{\small $\sim$onn}}

\vskip.5cm\noindent {\small Uriel G. Rothblum}\newline
\emph{Technion - Israel Institute of Technology, 32000 Haifa, Israel}

\emph{email: rothblum{\small @}ie.technion.ac.il}

\emph{http://ie.technion.ac.il/{\small rothblum.phtml}}

\end{document}